\documentclass[12pt,sort]{elsarticle}
\biboptions{sort&compress}
\usepackage[colorlinks=true, linkcolor=blue, citecolor=blue, urlcolor=blue]{hyperref}

\usepackage{amsmath}
\usepackage{amssymb}

\usepackage{algorithm}
\usepackage{algpseudocode}
\DeclareMathOperator{\st}{s.t.}

\newcommand{\OA}{\mathcal{O}}
\newcommand{\IA}{\mathcal{I}}
\newcommand{\inR}{\in\mathbb{R}}

\usepackage{tabularx}
\usepackage{threeparttable} %

\usepackage{booktabs}
\usepackage{subcaption}
\usepackage{graphicx}

\usepackage{xcolor}
\definecolor{orange}{rgb}{1,0.5,0}

\newcolumntype{L}{>{\raggedright\arraybackslash}X} %

\journal{Computers \& Chemical Engineering}

\begin{document}

\begin{frontmatter}

\title{ORACLE: A rigorous metric and method to explore all near-optimal designs for energy systems
}

\author[inst1]{Evren M. Turan}
\author[inst1]{Stefano Moret}
\author[inst1]{André Bardow}
\ead{abardow@ethz.ch}
\affiliation[inst1]{organization={Energy \& Process Systems Engineering},%
            addressline={ETH Zürich}, 
            city={Zürich},
            postcode={8092}, 
            country={Switzerland}}

\begin{abstract}

Optimization models are fundamental tools for providing quantitative insights to decision-makers. However, models, objectives, and constraints do not capture all real-world factors accurately. Thus, instead of the single optimal solution, real-world stakeholders are often interested in the near-optimal space -- solutions that lie within a specified margin of the optimal objective value. Solutions in the near-optimal space can then be assessed regarding desirable non-modeled or qualitative aspects. The near-optimal space is usually explored by so-called Modelling to Generate Alternatives (MGA) methods. However, current MGA approaches mainly employ heuristics, which do not measure or guarantee convergence. 

We propose a method called ORACLE, which guarantees generation and exploration on the \emph{entire near-optimal} space by exploiting convexity. ORACLE iteratively approximates the near-optimal space by introducing a metric that both measures convergence and suggests exploration directions. Once the approximations are refined to a desired tolerance, any near-optimal designs can be generated with negligible computational effort. 

We compare our approach with existing methods on a sector-coupled energy system model of Switzerland. ORACLE is the only method able to guarantee convergence within a desired tolerance.  Additionally, we show that heuristic MGA methods miss large areas of the near-optimal space, potentially skewing decision-making by leaving viable options for the energy transition off the table.
\end{abstract}


\end{frontmatter}

\section{Introduction}
The transition away from fossil resources necessitates massive shifts in the energy, chemicals, and fuels industries. However, decision-making in large-scale systems is inherently complex due to the intricate interdependencies among various decisions -- for example, in an energy system, the cost-effectiveness of low-carbon power-to-hydrogen technologies depends on not only their costs but also the availability of green electricity or alternative generation and storage technologies \citep{ajanovic2024future, wurbs2024important}. 

A common strategy to manage this complexity is to develop models that describe an abstracted version of a system, and run these models to provide insights to decision-makers. These models are often optimization-based techno-economic models that describe the technical, quantifiable requirements of the system (such as affordability, planning, and operation requirements) and find the best decisions to satisfy these concerns while minimising a specified objective, most commonly the total system cost. Thus, given a set of input parameters that make up a scenario, these models return the least-cost solution for that scenario. Energy system and integrated assessment models \citep{nordhaus1992optimal,manne1976eta,nordhaus1977economic,fishbone1981markal,PyPSAEur,limpens2019energyscope} are examples of such models that have been used since the 1970s to aid stakeholders in their decision-making. 

However, providing only a single cost-optimal solution is both limiting and naïve \citep{trutnevyte2016does}. Least-cost optimization is effective for identifying technically feasible scenarios while providing an estimate of the required costs. However, due to uncertainty of the future, simplifying assumptions, and idealisations, we can be confident that the cost-optimal solution is not the “real life” cost-optimal solution \citep{stigler1945cost,trutnevyte2016does}. Indeed, hard-to-model social and political factors can render the least-cost scenario untenable \citep{stigler1945cost}.  %

A particularly promising approach to handle these concern  %
is to identify \emph{near-optimal} solutions -- that is, solutions that are similar in cost but otherwise differ \citep{brill1979use,decarolis2011using, lau2024measuring}. The intuition is that, regardless of other stakeholder interests, excessively expensive decisions are not of interest. Once the near-optimal solution space has been defined, decisions can be analysed while taking un-modelled factors and trade-offs into account. %

To find near-optimal designs, researchers have proposed a class of methods called Modelling to Generate Alternatives (MGA) \citep{lau2024measuring,grochowicz2023intersecting,brill1979use,price2017modelling,decarolis2016modelling,pedersen2021modeling,decarolis2011using}.  MGA methods mainly employ heuristic strategies to explore the near-optimal space by repeated re-optimization of the system model. This exploration usually has two goals: (i) to identify maximally different solutions and (ii) to ensure a diversity of solutions within the entire near-optimal space so that the exploration does not yield a bias towards/against specific combinations. However, as system models are typically large and computationally expensive, MGA methods are usually run for a specified, low number of iterations. Consequently, recent reviews have highlighted that the returned solutions rarely have guarantees of coverage, diversity, or of being unbiasedly distributed within the near-optimal space \citep{lau2024measuring}. 

A recent variant, Modelling All Alternatives (MAA) \citep{pedersen2021modeling,grochowicz2023intersecting}, attempts to address these issues by monitoring convergence of the volume of the space; however, due to computational complexity, MAA is limited to very low-dimensional problems ($\lesssim 8$ exploratory variables). 

In this paper, we propose an approach, Optimization for the Rigorous Analysis of the Complete Landscape of Alternatives (ORACLE), to overcome these challenges by efficiently approximating the entire near-optimal space of a convex optimization model. Our method exploits convexity to iteratively refine outer and inner approximations of the near-optimal region. For this purpose, we introduce an interpretable convergence metric that measures how much of the space is left to explore. Once a desired approximation tolerance is achieved, we can sample near-optimal solutions at negligible computational cost. The sampling can be tailored to meet specific criteria, e.g., to be maximally diverse or uniformly distributed in the space. Such properties are difficult to guarantee when using the original system model due to computational restrictions. 

Furthermore, our proposed convergence metric quantifies how much of the space is left to explore in interpretable units, and can be applied to monitor the convergence of existing and future MGA methods. We compare the ORACLE approach with existing methods on a large-scale sector-coupled energy systems model of Switzerland and show that significant areas of the near-optimal space are not explored by existing MGA methods, resulting in important near-optimal decisions being overlooked. In contrast, with similar computational effort per iteration, our method guarantees finding the entire near-optimal space. Crucially, while presented here for an energy system model, our proposed approach applies to any convex optimization model.

The paper is organised as follows: in Section \ref{sec: background}, we briefly review existing MGA methods and establish the necessary mathematical background.  In Section \ref{sec: methodology}, we detail our proposed algorithm, which is then numerically demonstrated in Section 
\ref{sec: num-exps}, before concluding the paper with a brief discussion and outlook in Section \ref{sec: conclusion}.

\section{Background and formal definitions}\label{sec: background}

\subsection{Optimal and near-optimal designs}
The operation and design of an energy system can be written as an optimization problem, where the objective is to minimize the total system's cost subject to the constraints that define the system. To manage the computational complexity, a widely used approach is to set up the model as a large linear program (LP) \citep{fishbone1981markal,PyPSAEur,limpens2019energyscope} as these can be solved reliably and efficiently. An LP can always be written in the form:
\begin{subequations}\label{eq: system-model-LP}
\begin{align}
    v^* = \min_{x}\:\: & c^T x \\
    s.t.\:\: & Ax \leq b \\
             &  Fx = d
\end{align}   
\end{subequations}
where $v^*\inR$ is the optimum objective value, $x\inR^{n_x}$ is the decision vector, and $A\inR^{n_i\times n_x},\ b\inR^{n_i},\ F\inR^{n_e\times n_x},\ d\inR^{n_e},\ c\inR^{n_x}$ are model parameters that define the $n_i$ inequality constraints, $n_e$ equality constraints, and the system costs respectively. As convention, we use the asterisk to denote optimality, e.g. $x^*$ is an optimum solution with cost $v^*$. 
The optimization variable contains $n_x$ entries (typically $\gg10^5$) relating to both design/capacity decisions (e.g., how much hydro power to install) and operational decisions (e.g., how should the installed hydro power be used in each time step of the model).  These energy system models are typically very large, and depending on the granularity of the spatial and temporal resolution, can take minutes to hours to solve, e.g., see \citet{limpens2019energyscope,brochin2024harder}, and the references therein. 

Let $\epsilon\%\ge0$ be the maximum allowable percentage deviation from the optimum value. The near-optimal space, $\mathcal{X}_\epsilon$, is defined as all points that satisfy the original constraints, and are within the budget defined by $\epsilon$:
\begin{equation}\label{eq: define-full-nearopt-space}
    \mathcal{X}_\epsilon = \left\{x\: \middle|\: Ax \leq b,\: Fx = d,\:c^T x \leq v^* (1+ \epsilon) \right\}
\end{equation}

As the decision vector $x$ contains many entries, we typically do not want to explore the entire decision space; rather, we want to focus on the space of some key decisions. For example, for energy systems, we may only want to explore design decisions or the total annual use of a technology instead of hourly operating decisions. Let the decisions to explore be denoted by  $z\inR^{n_z}$, and the corresponding near-optimal region be $\mathcal{Z}_\epsilon$:
\begin{equation}\label{eq: set-noteation-Zeps}
    \mathcal{Z}_\epsilon = \left\{z\: \middle|\:Sx = z,\: Ax \leq b,\: Fx = d,\:c^T x \leq v^* (1+ \epsilon) \right\}
\end{equation}
where $S\inR^{n_z\times n_x}$ is a linear projection into the exploratory variable space.
In the rest of the paper, $\mathcal{Z}_\epsilon$ will be referred to as the near-optimal space. 

We assume that the upper and lower bounds of $z$ are available (these can be determined from the upper and lower bounds of $x$). Furthermore, we assume that (i) for any $z$ within these bounds the equality constraint $Fx=d$ can be satisfied and (ii) that $\mathcal{Z}_\epsilon$ has a non-empty interior, i.e. that the near-optimal space has an ``inside". These assumptions are not restrictive, as they will be naturally satisfied for a reasonably defined study. %

\subsection{Modelling to generate alternatives}
 As systems models are typically very large, and not always available in a manipulable form, it is not tractable to find  $\mathcal{Z}_\epsilon$  by projecting
 the constraints of \eqref{eq: set-noteation-Zeps} onto $\mathbb{R}^{n_z}$. Instead, researchers typically approximate $\mathcal{Z}_\epsilon$  by finding a point set of near-optimal designs, using Modelling to Generate Alternatives (MGA) approaches. MGA involves  solving  $m$ optimization problems to get a matrix of points $\mathcal{Z}_{points}=\left[z_1^T,\dots, z_m^T\right]$. As the vast majority of MGA papers are used with LPs, we do not consider methods that specialize on other problem classes.
 \begin{table}
    \small %
    \caption{A comparison of MGA algorithms in the literature with a selection of references describing their use and implementation.}
    \label{tab: mga-over-view}
    \centering
    \begin{threeparttable}
        \begin{tabularx}{\textwidth}{L c X X}
            Method & Parallelizable & Guarantees  & References\\
            \hline  \hline
            HSJ \newline \tiny{(Hop-Skip-Jump)}  & No  & - & \cite{brill1979use,brill1982modeling,chang1982use,decarolis2016modelling,decarolis2011using,lau2024measuring} \\
            \hline
            Random & Yes & Finds all vertices of $\mathcal{Z}_\epsilon$ with unlimited iterations & \cite{berntsen2017ensuring,lau2024measuring} \\
            \hline
            VMM \newline \tiny{(Variable min/max)} & Yes & Finds the minimum and maximum usage for each decision variable in finite iterations.\tnote{1} & \cite{neumann2021near,neumann2023broad,nacken2019integrated,schwaeppe2024finding,lau2024measuring}\\
            \hline
            ERG \newline \tiny{(Efficient Random Generation)} & Yes & - & \cite{chang1982efficient,trutnevyte2016does,sasse2019distributional}\\
            \hline
            SPORES \newline\tiny{(Spatially explicit practically optimal alternatives)} & Partially \tnote{2} & Variations that first perform VMM inherit its guarantees.  & \citep{lombardi2023redundant, lombardi2020policy,LCA_spores_2025}\\
            \hline
            MAA \newline\tiny{(Modelling All Alternatives)} & No  & Converges to $\mathcal{Z_\epsilon}$ with unlimited iterations. & \citep{grochowicz2023intersecting,pedersen2021modeling,lau2024measuring} \\ 
            \hline
            Manhatten MGA & No & Finds a point cloud and distance such that there is no additional point in $\mathcal{Z_\epsilon}$ more than this distance from the point cloud. & \citep{price2017modelling} \\
            \hline
            ORACLE \newline\tiny{(Optimization for the Rigorous Analysis of the Complete Landscape of Alternatives)} & Partially\tnote{3} & Converges to within a desired tolerance of $\mathcal{Z_\epsilon}$ in a finite number of iterations.  & This work. \\
            \hline
        \end{tabularx}
        \begin{tablenotes}
            \tiny
            \item[1] Some authors only explore either the maximum and/or minimum usage. If implemented to swap to Random after finding the min/max usages, then the guarantees of Random apply. 
            \item[2] Search directions are generated sequentially in the main part of the algorithm, but multiple sequences can be generated in parallel.
            \item[3] Calculation of the convergence metric has to be done sequentially, however many trial and feasible points can be generated in parallel. 
        \end{tablenotes}
    \end{threeparttable}
\end{table}
Most MGA algorithms solve the following problem at each iteration:
\begin{subequations}\label{eq:MGA}
\begin{align}
    \min_{x,z}\:\: & w_k^T z \\
    s.t.\:\: & Sx = z\\
            & Ax \leq b \\
             &  Fx = d\\
            & c^T x \leq v^* (1+ \epsilon)
\end{align}  
\end{subequations}
where $w_k\inR^{n_z}$ is the $k^{th}$ vector of weights chosen by the MGA algorithm and $W_{points}=\left[w_1^T,\dots, w_m^T\right]$ is the dataset of weights used. Most MGA algorithms only differ in how $W_{points}$ are chosen. 

As exact MGA implementations differ between papers, we describe prototypical MGA approaches below, with important aspects of these approaches summarised in Table \ref{tab: mga-over-view}.  Further details of MGA algorithms can be found in \ref{subsec: MGA-algs} or in the recent review of \citet{lau2024measuring}. 

    \begin{description}
        \item [Hop-Skip-Jump (HSJ):] \cite{brill1979use,brill1982modeling,chang1982use,decarolis2016modelling,decarolis2011using,lau2024measuring} \\
         HSJ is the earliest method known to us in the MGA literature. At each iteration, HSJ heuristically defines $w_k$ with the aim of finding a $z_k$ that is far from existing points in $\mathcal{Z}_{points}$. In HSJ, entries of  $w_k$ are chosen as the number of times the corresponding variable has a non-zero value in $\mathcal{Z}_{points}$. A variant, HSJ-rel, allocates entries as the sum of each entry in $\mathcal{Z}$ normalized by each variable's maximum attainable value. HSJ and HSJ-rel are inherently sequential, and hence MGA iterations cannot be computed in parallel. No convergence guarantees have been established for either variation. Additionally, prior authors have noted that HSJ tends to explore the space poorly \citep{lau2024measuring}. 
        
        \item [Random (vector):] \cite{berntsen2017ensuring,lau2024measuring}\\
        $w_k$ is chosen randomly from a uniform distribution, with entries between -1 and 1. Although this approach is very simple, it has two important properties: 
        \begin{enumerate}
            \item In high-dimensional space, any two random vectors are likely to be near-orthogonal \citep{blum2020foundations}. Thus, different parts of the space are naturally explored, and at first, each new weight is likely to be near-orthogonal to \emph{all} previous weights.
            \item As the number of iterations goes to infinity, this approach is guaranteed to explore every vertex of the near-optimal space. 
        \end{enumerate}

        \item [Variable Min/Max (VMM):]  \cite{neumann2021near,neumann2023broad,nacken2019integrated,schwaeppe2024finding,lau2024measuring}\\
        VMM \textit{systematically} assigns $w$'s with only one non-zero entry of either $1$ or $-1$ to find the maximal and/or minimal possible usage of each technology within the near-optimal budget. %
        Hence, at most $2n_z$ points will be found. After these points, the algorithm can either terminate or proceed to follow another MGA approach, e.g. Random, as in \citet{lau2024measuring}, if more points are desired. %
        
    \item [Efficient Random Generation (ERG):] \cite{chang1982efficient,trutnevyte2016does,sasse2019distributional}\\
        ERG \textit{randomly} assigns $w_k$ entries of -1, 0, or 1.
        Given infinite iterations, this method will recover the extreme solution bounds of VMM; however, unlike Random, it is not guaranteed to explore every vertex in the limit. %
        \item[Spatially explicit
practically optimal alternatives (SPORES)]\citep{lombardi2020policy,lombardi2023redundant, LCA_spores_2025}\\
            SPORES originates as a hybrid of HSJ-rel and VMM that was originally 
            proposed for spatially resolved models \citep{lombardi2020policy}.  Later variations have hybridised VMM with other weighting schemes \citep{lombardi2023redundant, LCA_spores_2025}. In these variations, the search direction is given by a weighted linear combination of the weight from VMM and the other schemes. Like HSJ-rel, these other schemes are sequential, meaning that an incumbent search direction is based on earlier search directions. If multiple sequences are generated, these sequences can be run in parallel. No convergence guarantees specific to SPORES have been established. If exploration is first performed using only VMM, and then the hybridisation is introduced (as in \citet{LCA_spores_2025}), SPORES inherits VMM's guarantees. %
        \item[Modelling all alternatives (MAA):] \cite{grochowicz2023intersecting,pedersen2021modeling,lau2024measuring,schricker2023gotta} \\
    Unlike the other MGA approaches, MAA is a region-based approach. MAA first generates an initial point cloud with VMM. After this phase, MAA iteratively forms the convex hull of $\mathcal{Z}_{points}$ and tries to expand this convex hull by selecting $w_k$ to be the normal of one of the facets of the hull. This iteration continues until a convergence criterion based on the volume of the convex hull is met. As discussed in Section \ref{subsec: region-based-mga}, this approach is limited to a low number of exploratory variables due to the computational effort of forming the convex hull. 
    \item[Manhattan MGA]\citep{price2017modelling}\\
        To the author's knowledge, Manhattan MGA is the only MGA algorithm (for continuous models) that does not use the standard MGA problem \eqref{eq:MGA}. Instead, the objective is to maximize the Manhattan distance between the new point and all previous points. This requires formulating an MILP with the energy system model and is thus computationally impractical for large models. 
    \end{description}

\subsection{Use of the convex hull in MGA}\label{subsec: region-based-mga}
The underlying idea of MAA is forming the convex hull of $\mathcal{Z}_{points}$ to explore the near-optimal region. The convex hull of a set of points is the smallest polytope that contains all the points.
Since the model is convex, the convex hull of $\mathcal{Z}_{points}$ forms the polytope, $\mathcal{I}$, which lies entirely within the near-optimal region, i.e. $\mathcal{I} \subseteq \mathcal{Z}_\epsilon$. Thus, $\mathcal{I}$ is an inner approximation, as all points in $\mathcal{I}$ are guaranteed to be near-optimal. Authors have proposed using the convex hull both to guide the exploration of $\mathcal{Z}_\epsilon$ (MAA), and as a post-processing step to generate further points and perform further analysis \citep{lau2024modelling}.

MAA requires forming the halfspace representation of the convex hull from the $m$ points in $\mathcal{Z}_{points}$. For more than a few dimensions, computing this representation is computationally intractable due to the curse of dimensionality: for state-of-the-art methods, the execution time grows $\mathcal{O}(m^{\lfloor n_z/2 \rfloor})$ \citep{barber2013qhull}. Thus, practically MAA-based methods are only applicable for $n_z\lesssim 8$, and even then can be computationally prohibitive when the number of points is large.

 In summary, most MGA methods explore $\mathcal{Z}_{points}$ by finding a fixed, pre-determined number of points without measuring convergence or coverage of the space. 
MAA and Manhattan can measure their coverage of the space\footnote{The original MAA method \citep{pedersen2021modeling} uses a heuristic based on relative improvement, but a later version of the method \citep{schricker2023gotta} introduces a valid convergence metric. See Section \ref{subsec: comparison-of-volume-metrics} for futher details.}; 
 however, both have computational limitations.  %
 Thus, no existing method can map the near-optimal space of a large-scale convex optimization model without restricting the exploration to only a few variables. 
 
\section{Methodology}\label{sec: methodology}
\subsection{The ORACLE algorithm}

To overcome the limitations discussed in Section \ref{subsec: region-based-mga}, we propose a method to map near-optimal spaces, called ORACLE: ``Optimization for the Rigorous Analysis of the Complete Landscape of alternatives in Energy systems"\footnote{The source code of ORACLE is available at: \url{place-holder-until-publication.}}. The underlying idea of ORACLE is to leverage convexity to \emph{adaptively} construct approximations of the near-optimal space while measuring the convergence of these approximations to within a desired tolerance. 

\begin{figure}
    \centering\includegraphics[width=.99\linewidth]{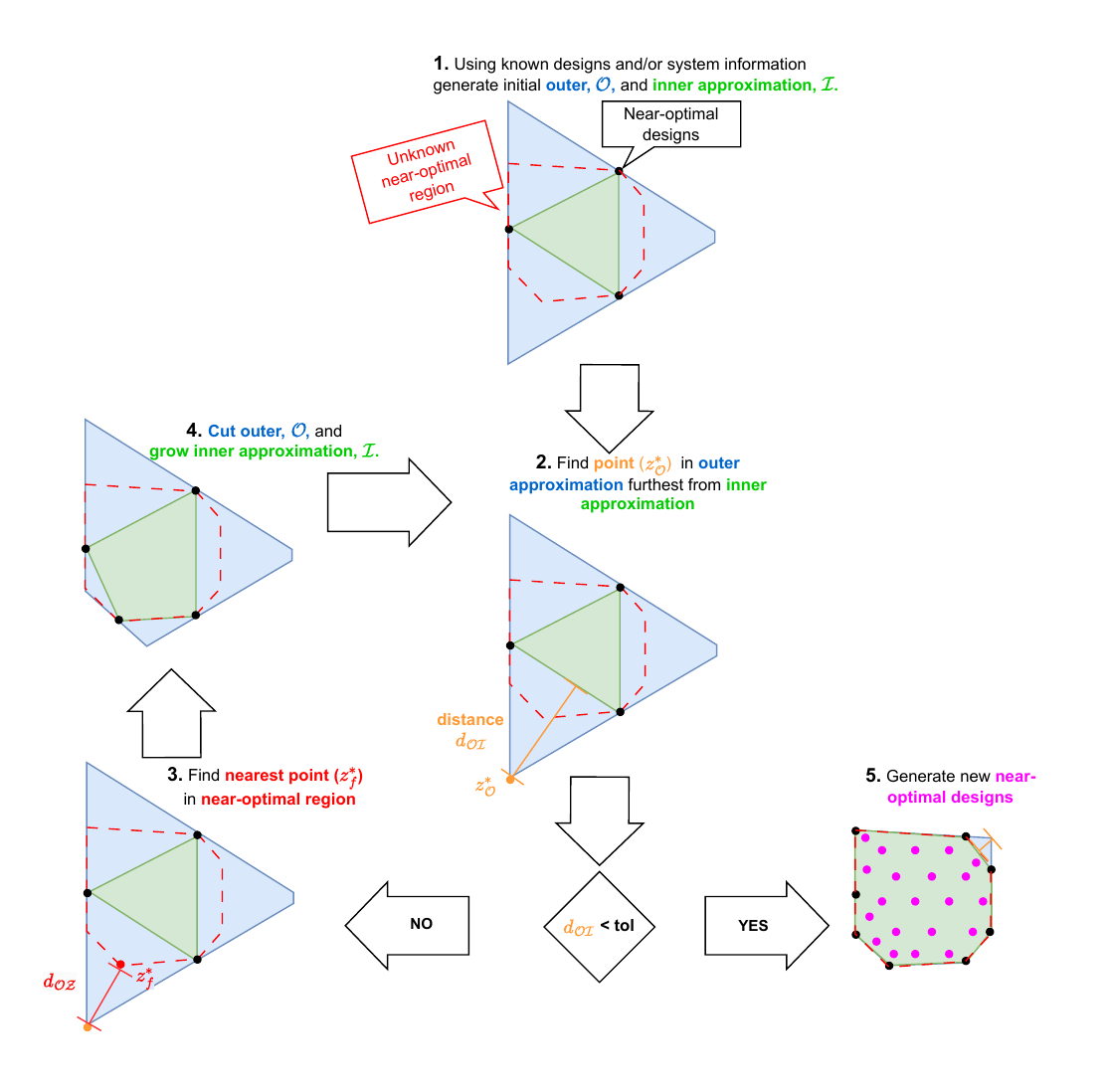}
    \caption{The key steps of the ORACLE algorithm. ORACLE is initialized by constructing an initial outer and inner approximation ($\OA$ and $\IA$) using known near-optimal points and system information. In Step 2, the main loop of ORACLE begins by finding the point in $\OA$ furthest from $\IA$. If the distance is below a desired tolerance, the algorithm terminates, and the resulting polytope(s) can be used to generate near-optimal designs at negligible computational effort. Otherwise, ORACLE locates the nearest feasible point to the outer point (Step 3), and uses this solution to grow the inner approximation, and (if possible) cut the outer approximation (Step 4), before going back to Step 2.}
    \label{fig:prop-alg}
\end{figure}

ORACLE, shown in Figure \ref{fig:prop-alg}, iteratively shrinks an outer approximation, $\OA$,  while iteratively expanding an inner approximation, $\IA$, until they converge to within some distance, $tol$. Importantly, this distance is the \emph{maximum possible} error incurred when generating near-optimal points, as the true near-optimal region is sandwiched between the two approximations. Thus after generation, the regions can be used, e.g. to generate many near-optimal points, while guaranteeing that the maximum error is below a desired tolerance.
Additionally, as shown in Section \ref{sec: num-exps}, this metric can also be used to measure the convergence of other MGA approaches.

The key steps of ORACLE are outlined below, with numbering corresponding to the steps shown in Figure~\ref{fig:prop-alg}. Further details are provided in the subsequent sections, while the pseudocode and potential variations are in \ref{subsec: oracle-pseudo}.

\begin{description}
\item [Step 1:] Given the cost-optimal design, and additional system information (e.g., other known designs, maximum capacities), using convexity, we form initial inner and outer approximations ($\IA$, $\OA$) of the entire near-optimal region.
\item [Step 2:] To monitor convergence, and to guide exploration of the near-optimal space, we find a trial point, $z_{\mathcal{O}}^*$, which is the furthest point in the outer approximation from the inner approximation. This distance ($d_{\OA\IA}$) is directly interpretable as the maximum error of the approximations.
If this distance is below a desired tolerance, we proceed to Step 5; otherwise, we proceed to Step 3. 
\item [Step 3:] We then explore in the direction of the trial point by finding the nearest point in the near-optimal region, $z_f^*$, to  $z_{\mathcal{O}}^*$. This requires solving the energy systems model with a distance-minimizing objective.
\item [Step 4:] We then grow the inner approximation, $\IA$, by including the new point, $z_f^*$, in its data. The outer approximation is shrunk, cutting off a portion of $\OA$ which contains $z_{\mathcal{O}}^*$. 
\item [Step 5:] The refinement of the approximations stops once the distance between the approximations ($d_{\OA\IA}$) falls below a desired tolerance.  At this point,  near-optimal designs can be generated or other analyses performed with negligible computational effort using only the inner and outer approximations.
\end{description}

\subsection{Mathematical formulation} 
\label{sec: methodology-region-based-approx}
This section provides in-depth details on the mathematical formulations and the steps of ORACLE.
\subsubsection{Preliminaries of inner and outer approximations}
 To form an inner approximation, $\IA\subseteq\mathcal{Z}_{\epsilon}$, we define the convex hull implicitly as the convex combination of all points in $\mathcal{Z}_{points}$:
\begin{equation}   \label{eq: define-IA}
\mathcal{I} = \left\{ z \:\middle|\: \lambda^T \mathcal{Z}_{points} = z, \:\: \sum_{i=1}^m \lambda_i = 1, \:\: \lambda \geq 0\right\}
\end{equation}
with $\lambda\inR^{m}$.
Thus, for a point $z$ to be in $\IA$, one must be able to find a $\lambda$ satisfying the constraints of \eqref{eq: define-IA}. 
Unlike explicitly constructing the convex hull, verifying whether a point satisfies \eqref{eq: define-IA} (or constraining a point to lie in  $\IA$) is an LP, which remains computationally tractable even for many high-dimensional points. While this formulation of the convex hull was used by \citet{lau2024modelling} for post-processing of MGA solutions, this is the first time it has been proposed to explore the space to the best of the authors' knowledge. 

In addition, we define an outer approximation of the near-optimal space, $\mathcal{Z}_{\epsilon}\subset\OA$, by defining a set of linear inequality constraints that are true for any  $z\in\mathcal{Z}_\epsilon$:
\begin{equation}
    \mathcal{O} = \left\{z \:|\: y_{O_i}^{T} (z-z_{O_i}) \leq 0, \:\: \forall i=1,\dots,n_{\mathcal{O}} \right\} \label{eq: def-OA}
\end{equation}
where $y_{O_i}\inR^{n_z}$ and $z_{O_i}\inR^{n_z}$ are parameters -- respectively the $i$th normal and point of a hyperplane. A simple way to form an outer approximation is to note that, due to convexity, the search direction $w_i$ used in the standard MGA problem \eqref{eq:MGA} is the normal of a supporting hyperplane of $\mathcal{Z}_\epsilon$ through $z_i$ \citep{boyd2004convex}. 

 In the following, we provide further details on the steps that make up ORACLE.
\subsubsection{Step 1: Defining the initial regions}

The inner approximation is always defined as the convex combination of known points, as in \eqref{eq: define-IA}. Initially, the cost-optimal solution will be known, so the inner approximation will always have at least one point. Any other known designs can be included in the definition. 

The initial outer approximation is defined by (linear) inequalities that are known to be satisfied for any point in the true near-optimal region. The simplest example of such constraints are the upper and lower bounds of $z$, (e.g. the first element may have bounds $0\le z_1 \le 5$). Another simple constraint that can be included is that an under-approximation of the system cost $\underline{v}$, should satisfy the near-optimality constraint\footnote{If $\underline{v}$ is an under-approximation, then any point in $\mathcal{Z}_\epsilon$ is guaranteed to satisfy \eqref{eq:under-approx-cost}. Thus, to be part of the over-approximation, $\underline{v}$ must be an under-approximation. }:
\begin{equation}\label{eq:under-approx-cost}
   \underline{c}z = \underline{v} \le v^*(1+\epsilon)
\end{equation}
where $\underline{c}\in\mathbb{R}^{n_z}$ is some parameter that under-approximates the contribution of $z$ to the total system cost. For example, often the choice $\underline{c}= Sc^T$ is valid as this approximates the total system cost as only the cost of using decisions $z$. Although this approximation is unlikely to be tight, it can be useful if an exploratory variable is unbounded (or has a very high bound), by ensuring that the outer approximation is bounded.
For simplicity, we assume that the initial outer approximation is compact (it is bounded and includes its boundary).

    \subsubsection{Step 2: Finding a trial point in the outer approximation} 
    The next step is to find a trial point in the outer approximation that is not in the inner approximation, and ideally outside the near-optimal region.  In ORACLE, we find the point in the outer approximation that is furthest from the inner approximation.
    The intuition is that either the trial point is in the near-optimal region, and hence the inner approximation should grow greatly (as this point was far from it), or the trial point is outside the near-optimal region, and the outer approximation should shrink. This refinement is done in the following step. 
    
    To generate this trial point, we maximise the minimum distance between $\IA$ and $\OA$:
\begin{equation}\label{eq: maxmin_dist_prob}
        d_{\mathcal{IO}} =\max \limits_{ z_O \in \mathcal{O}}\:\:\min_{z_I\in \mathcal{I}}\:\: \|z_{O} - z_{I}\|_p 
\end{equation}
where $p$ determines the norm used to measure distance. This problem returns $d_{\mathcal{IO}}$, the max-min distance between $\IA$ and $\OA$, the trial point, $z^*_O$. As the near-optimal region lies between the outer and inner approximation, $d_{\mathcal{IO}}$ is an upper bound on the distance to the near-optimal region. Thus, it is interpretable as the maximum error of using either $\mathcal{I}$ or $\mathcal{O}$ to approximate the near-optimal region $\mathcal{Z}_\epsilon$.

Problem \eqref{eq: maxmin_dist_prob} is a bilevel optimization problem and cannot be provided as is to most solvers. %
By picking a $p$ such that the lower level is convex, e.g. $1,\: 2,\: \text{or } \infty$\footnote{Also known as the maximum norm, as $\|z_{O} - z_{I}\|_\infty = \max \limits_{i=1,...,n_z} |z_{O_i} - z_{I_i}|$.}, problem \eqref{eq: maxmin_dist_prob} can be reformulated to yield a single-level problem.

We choose the infinity norm ($p=\infty$) as this measures the maximum error due to a single element of $z$, making the metric easy to interpret.
With this choice, problem \eqref{eq: maxmin_dist_prob} is a bilevel LP, which we reformulate into a single-level mixed integer linear program (MILP, see \ref{subsec: max-min reformulations} for details on the reformulation), which can be solved using off-the-shelf solvers.

    \subsubsection{Step 3: Finding the closest feasible point to the trial point}
    Given a trial point, we need to check if it is feasible. To do so, we use the system model to find the closest feasible point, using the same distance metric as in problem \eqref{eq: maxmin_dist_prob}. %
    Selecting $p=\infty$, yields the problem:
       \begin{subequations}\label{eq: closest-feasible-inf} \begin{align} 
        \min_{\delta,z_f,x}\:\: & \|\delta\|_\infty \\
            &  z_{\mathcal{O}}^*-z_f = \delta \label{eq: define-delta-for-duals}\\
                & Sx = z_f\label{eq: define-delta-for-duals-2}\\
                & Ax \leq b \\
                 &  Fx = d\\
                & c^T x \leq v^* (1+ \epsilon)
    \end{align}     \end{subequations}
where $\delta$ is the difference between between $z_f$ and $z_{\mathcal{O}}^*$, and $z_f^*$ is the optimal solution. Reformulating the infinity norm as linear inequalities results in an LP of similar size to solving the original systems model \eqref{eq: system-model-LP} or the standard MGA problem \eqref{eq:MGA}.

    If the trial point $z_{\mathcal{O}}^*$ is infeasible, we can both refine the inner approximation (using the nearest feasible point) and generate a cut that will be a supporting hyperplane of the true near-optimal region (as $z_{f}^*$ lies on the boundary). If we directly checked the feasibility of the trial point  $z_{\mathcal{O}}^*$\footnote{i.e. by using the constraint $Sx = z_{\mathcal{O}}^*$.}, then (i) the inner approximation would only be refined if $z_{\mathcal{O}}^*$ proved feasible, and (ii) if $z_{\mathcal{O}}^*$ were infeasible the resulting cut might be positioned far from the boundary of the near-optimal region, limiting its effectiveness.

    \subsubsection{Step 4: Refining the outer and inner approximations}

    At each iteration, we can grow the inner approximation by including the closest feasible point, $z_f^*$, in the set of points defining the polytope \eqref{eq: define-IA}. In contrast, if the trial point is feasible, $z_f^*=z_{\mathcal{O}}^*$, then we cannot refine the outer approximation at this iteration as this point ``confirms" the current outer approximation.  When the trial point is infeasible,  the outer approximation can be refined by adding a valid inequality to equation \eqref{eq: def-OA}.%
    
    If $z_f^*$ is a vertex of the constrained region, then infinitely many supporting hyperplanes exist at that point. Our goal is to find a \emph{strictly separating hyperplane}, i.e. one that satisfies $y_{O_i}^Tz_{\mathcal{O}}^*>y_{O_i}^Tz_f^*$. This requirement ensures that whenever an infeasible point is found, the trial point, and thus part of the outer-approximation $\OA$, is cut off by introducing the new constraint into $\OA$. This property guarantees convergence of $\OA$, as the initial $\OA$ has finite volume (being compact) and each iteration removes a nonzero volume from it. %
    
    A strictly separating hyperplane is given by the dual variables of equation \eqref{eq: define-delta-for-duals} -- see \ref{subsec: hyperplane-from-dual} for further details. We additionally cut $\OA$ by introducing a linear under-approximation of the total system cost at $z_f^*$. This cut requires an additional minimization of the total cost with the constraint $Sx=z_f^*$, but can accelerate convergence to the near-optimal region (see \ref{subsec: double-cut-OA}). %

 \subsubsection{Step 5: Sampling from the approximations}
 Once the approximations have converged to a desired tolerance, one can sample from either the inner ($\IA$) or outer approximations  ($\OA$) with considerably less computational effort than it takes to solve the original energy system model due to the dramatic decrease in problem size ($n_z\ll n_x$). The sampling can be performed very flexibly, and allows for the potential use of (i) any MGA method with the $\IA$ or $\OA$ as constraints instead of the original model as constraints, (ii) Markov chain Monte Carlo methods such as Hit-and-Run and its later variations  \citep{berbee1987hit,smith1984efficient,kaufman1998direction} to find points \emph{uniformly distributed} in the near-optimal region or (iii) any other optimization problem to find points with certain characteristics. For illustration, we find \textit{maximally distant points} in Section \ref{subsec: left-out}.

\subsection{Related work and limitations}

To the authors' knowledge, the idea of an outer approximation has only appeared once in the MGA literature \citep{schricker2023gotta}, where the volume of $\OA$ was used to normalise the volume of the convex hull to monitor convergence of MAA. Similarly to finding the convex hull, calculating the volume rapidly becomes computationally intractable \citep{barber2013qhull}, and is only practical for $n_z\lesssim 8$.  %

In contrast, ORACLE requires an MILP to find the trial point in the outer approximation in problem \eqref{eq: maxmin_dist_prob}. 
The computational cost of solving an MILP is often strongly determined by the number of integer variables in the problem. The proposed reformulation uses $2n_z + m$ integer variables, where $m$ is the number of points used to define the inner approximation. Thus, for problems that require many points for convergence, the MILP may become computationally costly. Methods for reformulating and solving bilevel LPs are an active area of research \citep{bard2013practical, pineda2018efficiently, kleinert2021computing}. We have not compared the computational effort of these different formulations, nor in the potential use of heuristic methods to solve the MILP resulting from problem \eqref{eq: maxmin_dist_prob}.

ORACLE adaptively approximates a convex polytope. This geometric problem has been studied in other contexts, including the compression of large systems of linear inequalities, determining reachable sets, multi-objective optimization and scheduling \citep{kamenev1992-adaptive-approximation-polytope, mcclure1975polygonal,lotov1989generalized,efremov2009properties, sung2007attainable}. The major differences of ORACLE to these works are (i)  the computation and formulation of the proposed distance metric to find a trial point, (ii) the target polytope is indirectly defined by a very large optimization problem where the problem data is challenging to manipulate, and (iii) the application to large-scale energy systems. We note that the approach proposed in \citet{sung2007attainable} for finding the attainable region of a scheduling problem combines elements of ORACLE and MAA, in that a variant of problem \eqref{eq: maxmin_dist_prob} is solved with $p=2$, and the search direction restricted to be a normal of the facets of the inner approximation. This requires forming the convex hull and swapping between half-space and vertex representations, as in MAA, which rapidly becomes computationally intractable \citep{barber2013qhull}.

 \section{Numerical experiments}\label{sec: num-exps}
 \subsection{Model and experimental set up}
We use the EnergyScope TD model of Switzerland \citep{limpens2019energyscope} to demonstrate the performance of the ORACLE algorithm and to benchmark it against several existing MGA approaches: Random, VMM, HSJ, SPORES, and ERG. We exclude MAA and Manhattan MGA from the comparison due to their computational limitations. Additionally, we evaluate our proposed metric against the volume-based metric found in the MGA literature \citep{lau2024measuring,schricker2023gotta, pedersen2021modeling}. All optimization problems are implemented in AMPL \citep{fourer1990ampl} and solved with Gurobi version 12 \citep{gurobi}. When solving LPs, we use Gurobi's homogeneous barrier algorithm method without crossover and set the convergence tolerance to 1e-6. For MILPs, we use a relative and absolute MIP gap of 0.1  and 0.05, 32 threads, and impose a time limit of 600 seconds. Volumes are calculated with Quickhull \citep{barber2013qhull} using the ``QJ" setting.
 
 EnergyScope TD is a sector-coupled energy system model that optimizes a greenfield representation of Switzerland in 2050. The model is spatially aggregated into a single node and resolved hourly using twelve typical days.  We selected EnergyScope TD for this study as it is a national-scale sector-coupled model that is computationally efficient (evaluations of $\le60$ seconds), making it suitable for extensive comparisons across MGA methods.
 
 In our numerical results, we consider the following six installed capacities as exploratory variables: photovoltaics (PV), (onshore) wind, nuclear, natural gas combined cycle (gas), with and without carbon capture and storage (CCS), and coal power plants (with CCS)\footnote{Although coal power plants do not appear in a cost-optimal system, they can be built in the near-optimal space.}. We emphasize that this list does not contain all technologies in the model.  A small list of technologies is necessary for the computation of volumes in the method comparison. All exploratory variables are measured in GW, and all MGA methods are run until they converge within 0.1 GW, calculated by equation \eqref{eq: maxmin_dist_prob}, or reach 200 iterations. A deviation of 0.1 GW is typically insignificant in the context of strategic decision-making for a national energy system; therefore, this tolerance is relatively stringent, ensuring that the near-optimal space has been rigorously explored.

 This section is arranged as follows: we compare ORACLE with other MGA methods, using first our proposed metric, equation \eqref{eq: maxmin_dist_prob}, and then by volume. %
 We then highlight the benefit of our method by using the converged approximation to identify designs that are missed by the other MGA methods. 

\subsection{Comparison of MGA convergence}

\begin{figure}
    \centering\includegraphics[width=.99\linewidth]{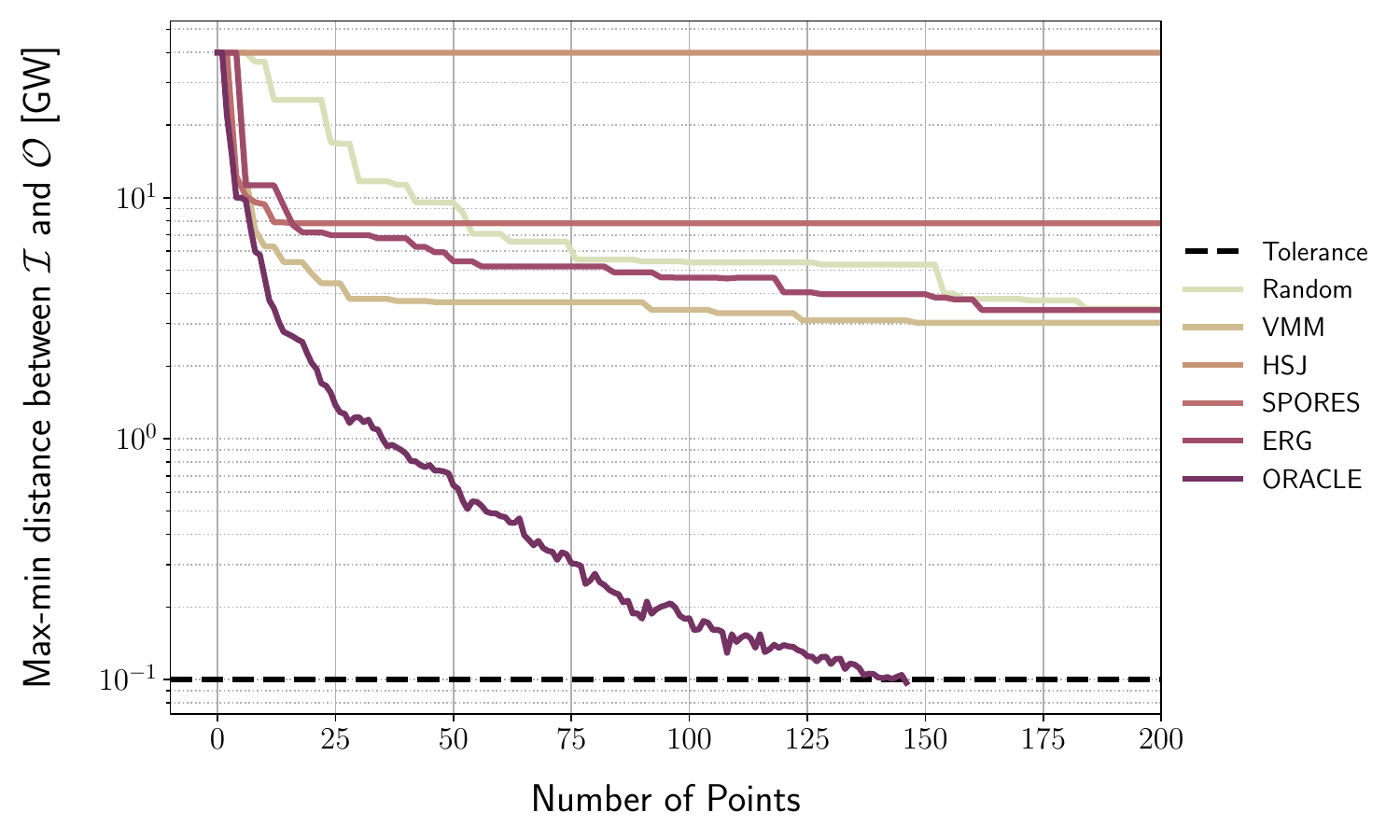}
    \caption{Distance between inner and outer approximations of the near-optimal space for state-of-the-art MGA methods and the ORACLE algorithm. ORACLE's non-monotonic convergence is a numerical artifact, as problem \eqref{eq: maxmin_dist_prob} is solved with relative and absolute termination criteria of 10\% and 0.05 respectively.}
    \label{fig:distance-oracle-small-tech-list}
\end{figure}
Although most MGA algorithms do not construct outer and inner approximations, they generate data that we use to construct these approximations (see Section \ref{sec: methodology-region-based-approx}). We initialize the approximations of all MGA algorithms using Step 1 of ORACLE. %

The methods (Figure \ref{fig:distance-oracle-small-tech-list}) show a stark contrast. ORACLE consistently outperforms all others and converges rapidly. ORACLE reaches a 1 GW tolerance in just 30 iterations, and converges to the desired tolerance (0.1 GW) within 144 iterations. None of the other methods even reach the 1 GW threshold within 200 iterations. The convergence of ORACLE is ``noisy" compared to the other methods due to the numerical tolerances used when solving problem \eqref{eq: maxmin_dist_prob}, which becomes visible when the distance becomes small.

HSJ shows minimal improvement, maintaining a nearly constant distance throughout. Similarly, SPORES initially improves rapidly before reaching a plateau. The maximum distance will only change\footnote{We emphasize that whenever the maximum distance remains constant, this does not mean that the size of the approximations necessarily remains constant (they could be growing or shrinking in other directions).} if a search direction ($w_k$ in problem \eqref{eq:MGA}) is picked that coincides with the direction of maximum distance. Both HSJ and SPORES generate a sequence of search directions that are autocorrelated (the search direction at iteration $k$ depends on previous weights). This autocorrelation may cause HSJ and SPORES to explore only a restricted part of the near-optimal space. Another group (Random, VMM, ERG) exhibits sporadic progress, reducing the maximum distance from 40 GW to around 4.5 GW within 200 iterations. The initial, rapid improvement of VMM is due to the systematic optimization for the ``real" range of each decision variable. After this phase, VMM is equivalent to the Random method, which has a chance to select the direction of maximum separation. As the ERG and Random methods are also similar to each other, it is unsurprising that these methods have similar trajectories.

\subsection{Comparison of MGA convergence by volume}\label{subsec: comparison-of-volume-metrics}
\begin{figure}
    \centering\includegraphics[width=.99\linewidth]{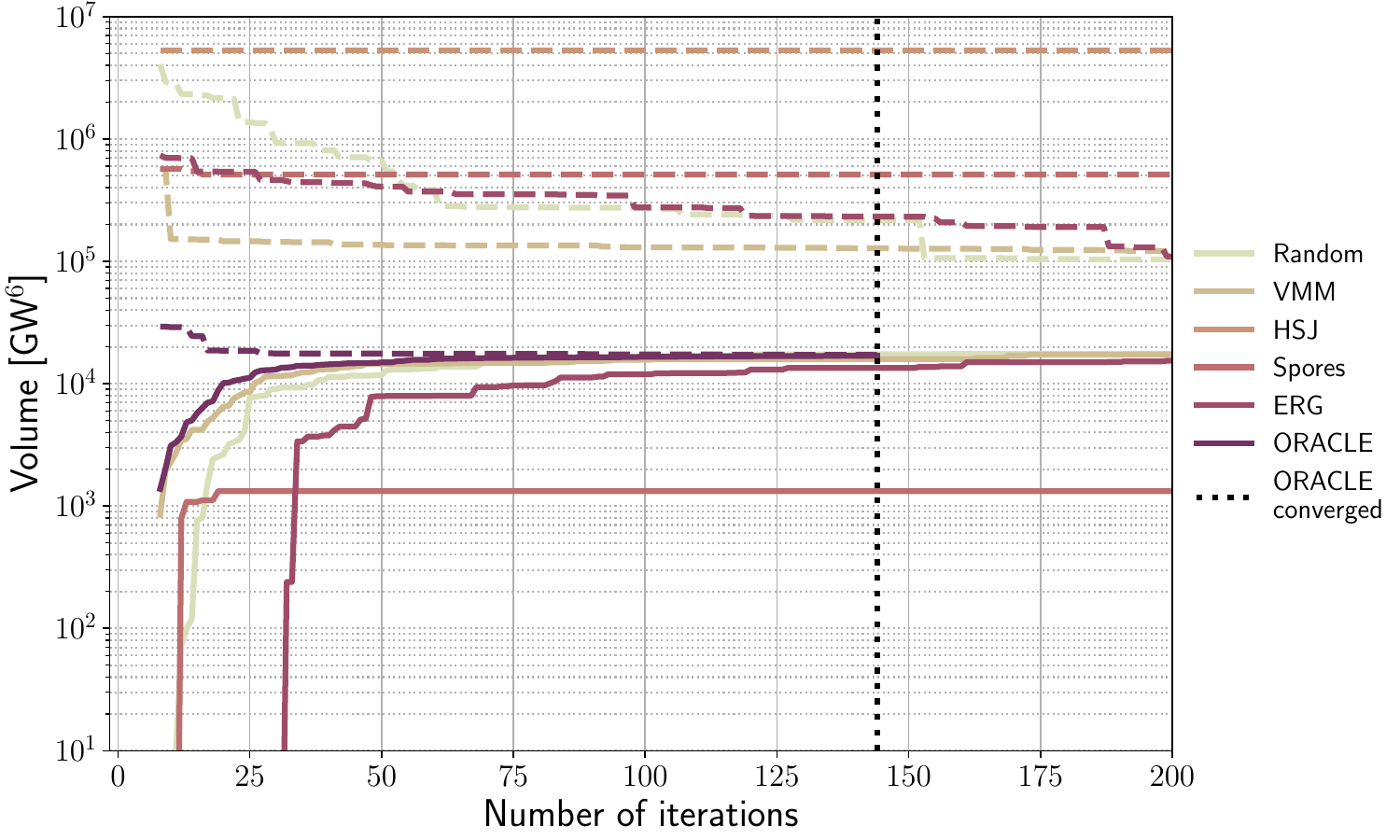}
    \caption{Volumes of the inner (solid lines) and outer (dashed lines) approximations of state-of-the-art MGA methods and of the ORACLE algorithm. ORACLE terminates after 144 iterations, having converged to within 0.1 GW (see Figure \ref{fig:distance-oracle-small-tech-list}).}
    \label{fig:compare-volumes-small-tech-list}
\end{figure}

Volume has been suggested as a convergence metric in the recent MGA literature \citep{pedersen2021modeling,lau2024measuring,schricker2023gotta}. When using the difference between the volumes of the inner and outer approximations as a metric, the groupings of the algorithms (Figure \ref{fig:compare-volumes-small-tech-list}) are the same as when the proposed metric (Figure \ref{fig:distance-oracle-small-tech-list}) is used: ORACLE performs the best followed by Random, VMM, ERG, and Spores with HSJ performing the worst. In fact, the volumes of the inner approximations of HSJ are not shown as they are very small in comparison. %
Small volumes have been reported for HSJ have been reported in earlier work \citep{lau2024measuring}.  

Interestingly, the volumes of the inner approximations of Random, VMM, and ERG are similar to the converged volume of ORACLE. However, as the volumes of their outer approximations are much larger, they are unable to show that the inner approximations appear to have converged by volume.  This finding implies that the poor performance of Random, VMM, and ERG in Figure  \ref{fig:distance-oracle-small-tech-list} is because their outer approximations are not tight. Despite the volumes being similar, the distance between the inner approximation of Random, VMM, and ERG with the converged outer approximation of ORACLE remains above the desired tolerance by a significant margin. Thus, although the volumes of the inner approximations are similar, the inner approximations have not converged to the desired tolerance (see \ref{subsec: compare-distance-true}).

Earlier work \citep{pedersen2021modeling,lau2024measuring} only used the change in the volume of the inner approximation as a heuristic to monitor convergence.  Figure \ref{fig:compare-volumes-small-tech-list} clearly illustrates the difficulty of using the rate of change, as it becomes small for ORACLE, Random,  VMM and ERG after approximately 50 iterations despite the inner approximations having different sizes. Additionally, as the volumes of the inner approximations of HSJ and SPORES are smaller than those of the other algorithms, their volumes change by less, which, by this metric, would be seen as convergence. Thus, the rate of change of volume does not appear suitable for measuring convergence. 
\citet{schricker2023gotta} proposed monitoring the ratio of the volumes of the inner and outer approximations, which is plotted in \ref{subsec: compare-volume-ratio}. This metric gives the same final ranking as our proposed distance metric; however,  due to computational limitations, volume can only be calculated for a small number of exploratory variables\footnote{In \citet{lau2024measuring}, it is suggested to approximate the volume by calculating the sum of the projected volumes of each pair of elements in $z$. This estimate is not guaranteed to increase whenever the actual volume increases, and so exacerbates the issues with using the rate of change.}.

\subsection{What designs are left out by the other MGA methods?}\label{subsec: left-out}
The proposed ORACLE method generates near-optimal regions that enable analyses that would be computationally infeasible using the original system model (defined by problem \eqref{eq: system-model-LP}). To demonstrate this capability, we formulate an optimization problem to answer the question: “What is the most different design that another MGA method did not find?”\footnote{As a further example of flexibility, note that a variation of problem \eqref{eq: most-distant-design} can be used to find the most diverse set of designs within the near-optimal region.}. This analysis is done by solving the following optimization problem to find the furthest point in the (approximated) near-optimal region to the point cloud generated by another MGA method:
\begin{subequations}\label{eq: most-distant-design}
    \begin{align}
        \max_{z_\mathcal{O}, \delta}\:\:& \delta \\
        &\delta \le d(z_\mathcal{O}-z_{MGA,i})\qquad \forall i=0,\dots,m\\
        &z_\mathcal{O} \in \mathcal{O}
    \end{align}
\end{subequations}
where $d$ is the desired distance function, $m$ is the number of points found by the other MGA methods, and $\delta$ is the maximum distance. In this section, we select the Manhattan distance, $d=\left\lVert\cdot\right\lVert_1$. Compared to optimizing the two and infinity norm, the Manhattan distance tends to find solutions that differ across multiple dimensions rather than finding solutions that differ greatly in only one dimension. 

\begin{figure}
    \centering
    \begin{subfigure}[b]{0.8\textwidth}
        \includegraphics[width=0.99\textwidth]{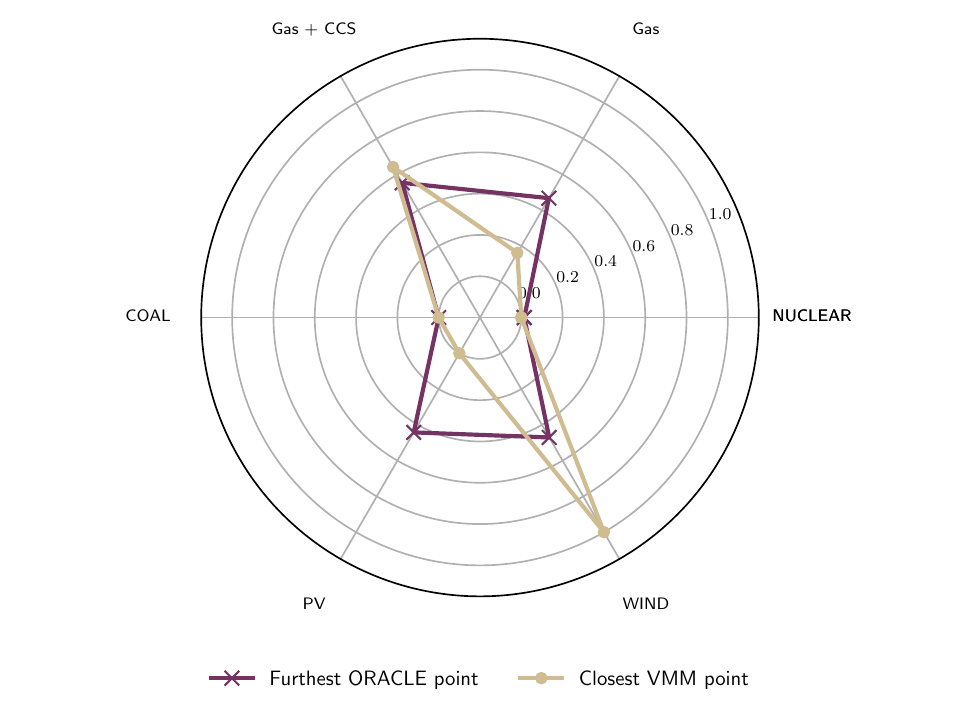}
        \caption{Comparison with VMM}
        \label{fig:radar-left-out-sub1}
    \end{subfigure}
    \vfill
    \begin{subfigure}[b]{0.8\textwidth}
        \includegraphics[width=0.99\textwidth]{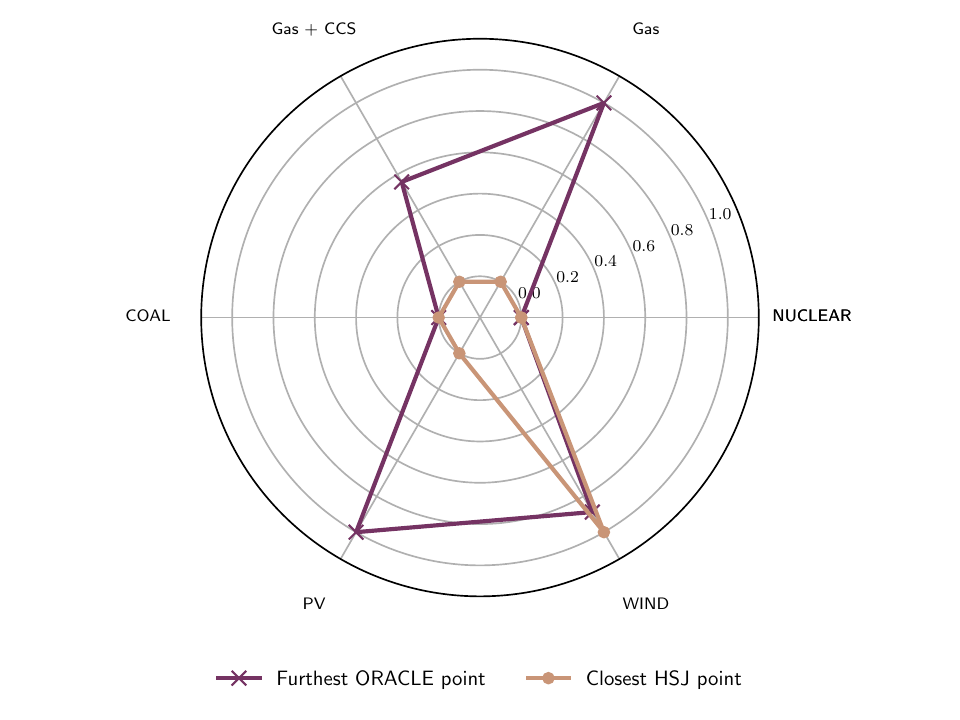}
        \caption{Comparison with HSJ}
        \label{fig:radar-left-out-sub2}
    \end{subfigure}
    \caption{Identification of strategies not found by the literature MGA methods. Each radial plot identifies a design found by ORACLE and the closest design found by the respective MGA method. The installed capacity of each technology is normalised against its maximum potential capacity.}
    \label{fig:radar-left-out}
\end{figure}

We select VMM and HSJ (run for 200 iterations) as representative examples from the ``groups" of Figures \ref{fig:distance-oracle-small-tech-list} and \ref{fig:compare-volumes-small-tech-list-ratio}, and identify their potentially incorrect conclusions due to missed designs (Figure \ref{fig:radar-left-out}) . For example, Figure \ref{fig:radar-left-out} (a) shows a near-optimal design identified by ORACLE that employs  moderate photovoltaics, wind, gas, and gas with CCS. However, the closest point found by VMM is one with high wind and no photovoltaics\footnote{The two designs have a significant difference in installed technologies. As not all electricity generation technologies are explored and energy imports are allowed, the HSJ design can use these to supply the energy demands of the system.}. Thus, one could incorrectly conclude that moderate deployment of photovoltaics, wind and gas is not viable. Figure \ref{fig:radar-left-out} (b) shows a more dramatic result; due to the poor convergence of HSJ, the closest HSJ point to the identified ORACLE design of high gas, wind and PV is very far away. In this case, one could incorrectly conclude that high deployment of wind implies that gas and photovoltaics should not be deployed, while their co-deployment is actually within the near-optimal space. 

This finding highlights that the lack of convergence metrics and exploration guarantees of MGA methods is not solely a theoretical concern, but a practical concern, as it can significantly impact the information provided to stakeholders. Additionally, as the proposed convergence metric, equation \ref{eq: maxmin_dist_prob}, directly estimates the maximum error in the results, its usefulness is immediately apparent: even if exploration is performed using a different MGA method, the proposed metric can be used to certify whether the results can be trusted. For example, as ORACLE has converged within 0.1 GW, the near-optimal space has essentially been fully explored.

\section{Conclusion}\label{sec: conclusion}

To avoid the bias of least-cost solutions \citep{stigler1945cost}, numerous researchers have proposed to instead generate \textit{near-optimal solutions} \citep{trutnevyte2016does, brill1979use, price2017modelling, decarolis2011using, schwaeppe2024finding}. Once generated, the near-optimal solutions can be analysed based on stakeholder interests, while taking un-modelled decision factors into account. Currently, near-optimal solutions are generated by a class of methods called Modelling to Generate Alternatives (MGA). These methods employ heuristics to generate diverse solutions in the near-optimal space. However, the majority of MGA methods do not measure their coverage of the space, have few or no performance guarantees, and normally run for a pre-specified number of iterations before terminating. Recent variants \citep{pedersen2021modeling, price2017modelling} have attempted to address these concerns, but have computational limitations that restrict their usage to either small models or a low number of exploratory variables.

We have proposed a convergence metric that quantifies how much of the space is left to explore. Importantly, this metric can be applied to existing MGA methods, allowing for easy estimation of their convergence. Furthermore, we have developed an approach called ORACLE that leverages convexity alongside the convergence metric to adaptively approximate the entire near-optimal space to within a desired tolerance. Once converged, near-optimal designs can be generated with negligible computational cost. In this step, the design generation can be tailored to have desired characteristics: e.g., to be most diverse, uniformly distributed, and so on. %

We compare ORACLE with five other MGA algorithms to explore the 10\% near-optimal region of a greenfield, sector-coupled energy model of Switzerland. ORACLE converges within 144 iterations, while the other MGA methods are unable to guarantee their convergence within 200 iterations (Figure \ref{fig:distance-oracle-small-tech-list}). We demonstrate the impact of this non-convergence by finding very different near-optimal designs that are missed by the MGA methods (Figure \ref{fig:radar-left-out}).  Additionally, we compare our proposed metric with the use of volume, as suggested in the MGA literature. While the rankings of the methods are the same, our proposed method can be calculated for larger dimensions than volume, while additionally providing a direct estimate of the error of the corresponding MGA method in convenient units.

The computational time per iteration of ORACLE is more expensive than that of other MGA methods, as it involves not only evaluating the system model, but also performing another optimization to evaluate the convergence metric. System models with high spatial and temporal resolution are computationally expensive to evaluate; hence, this step is likely to dominate the time per iteration, with the additional time per iteration of ORACLE being negligible in comparison (see Section \ref{subsec: comp-time}). To compensate for the increase in the computational time of using a larger model or exploring more variables, strategies can be employed such as partial parallelization and hybridization of ORACLE with other algorithms (for further details, see \ref{subsec: oracle-pseudo}). 

In conclusion, while the MGA literature has been rapidly growing \citep{trutnevyte2016does,lau2024measuring}, insufficient attention has been paid to the quality of exploration of MGA methods. Using the proposed metric, problem \eqref{eq: maxmin_dist_prob}, we can rigorously measure the potential error of using a particular method, thus allowing for MGA methods to be rigorously compared against each other. Thus, even if exploration is performed using a different MGA method, we recommend using the proposed metric to measure the exploration quality. 

\section*{Declaration of competing interest}

The authors declare the following financial interests and personal relationships that could be considered as potential competing interests: AB has served on review committees for research and development at ExxonMobil and TotalEnergies, companies active in both oil and gas and chemical production. AB holds ownership interests in firms that provide services to industry, some of which may operate in the chemical industry. In particular, AB has ownership interests in Carbon Minds, a company that supplies life cycle assessment (LCA) databases used, among others, in this work to validate the CRYSTAL model. The remaining other authors declare no known competing financial interests or personal relationships that could have influenced the work reported in this paper.

\section*{Acknowledgements}
TotalEnergies OneTech is acknowledged for partial funding through project CT0072000.
S.M. acknowledges support from the Swiss National Science Foundation under Grant No. PZ00P2\_202117. We thank Francesco Lombardi for his time and discussions regarding the implementation of SPORES.

\bibliographystyle{elsarticle-num-names} 
 \bibliography{bib}
 \newpage
\appendix
\section{Appendix}

\subsection{MGA algorithms}\label{subsec: MGA-algs}
In this section, we briefly review the implementation of the MGA algorithms used in the comparisons of the paper. Note that as exact implementations do differ between papers, we have implemented typical versions of these algorithms. These algorithms all generate a matrix of points $\mathcal{Z}_{points}=[z_1^T,\dots,z_m^T]$, where each $z\inR^{n_z}$ is a vector of the exploratory variables. At each iteration, $k$, the weighting vector $w_k\inR^{n_z}$ is specified, and problem \eqref{eq:MGA} is solved to yield solution $z_k$. The algorithms differ on how $w_k$ is chosen at each iteration.

    \begin{description}
        \item [Hop-Skip-Jump (HSJ):]\cite{brill1979use,brill1982modeling,chang1982use,decarolis2016modelling,decarolis2011using,lau2024measuring} \\ HSJ assigns weights based on how often a technology is used in the previous iterations. For this purpose, an element-wise indicator function, $\mathbf{1}_e$, equation \eqref{eq: elem-wise-indicator}  is used to define weights. This function returns a vector with elements that are one if the input element is non-zero, and zero otherwise. At each iteration weights are determined by: %
                \begin{subequations}
            \begin{align}
                w_1 &=\mathbf{1}_e(z_0^*)\\
                w_k &= w_{k-1} + \mathbf{1}_e(z^*_{k-1}), \qquad k=2,\dots,m\\
                \mathbf{1}_e(z)&=\begin{bmatrix}
1\:\: if\:\: z^{(1)} \not=0,\:\: else\:\: 0\\
\vdots\\
1\:\: if\:\: z^{(n_z)} \not=0,\:\: else\:\: 0\label{eq: elem-wise-indicator}
\end{bmatrix}
            \end{align}
        \end{subequations}
            HSJ-rel is a variation of HSJ where weights are updated by dividing the solution vector by a vector of the maximal values of each element:
                                \begin{subequations}
            \begin{align}
                w_1 &=\frac{z_0^*}{z^{max}}\\
                w_k &= w_{k-1} + \frac{z_{k-1}^*}{z^{max}}, \qquad k=2,\dots,m
            \end{align}
                    \end{subequations}
            where $z^{max}$ is a vector of the maximum values of the exploratory variables. 
            \item [Random (vector) MGA] \cite{berntsen2017ensuring,lau2024measuring}\\
            At each iteration, a vector of weights is generated with elements randomly between $-1$ and $1$.
            \item [Variable Min/Max (VMM):]  \cite{neumann2021near,neumann2023broad,nacken2019integrated,schwaeppe2024finding,lau2024measuring}\\
             VMM first determines the maximum and/or minimum possible usage of each technology within the near-optimal budget. A naive implementation requires solving at least $2n_z$ MGA problems; however, as the upper and lower bounds of the variables are normally known, the minimization or maximization of a variable can be skipped if it achieves its lower or upper bound in an earlier iteration. After this phase, we assign weights randomly, as in Random. Note that some other implementations of VMM stop after the initial phase. 
    
         \item [Efficient Random Generation (ERG):] \cite{chang1982efficient,trutnevyte2016does,sasse2019distributional}\\
         At each iteration, a random number of technologies is chosen to be randomly assigned a weight of $-1$ or $1$. All other technologies are assigned a weight of zero. 
        \item[Spatially explicit
practically optimal alternatives (SPORES)] \citep{lombardi2023redundant, lombardi2020policy,LCA_spores_2025}\\
        SPORES was proposed for spatially resolved models, and so we have adjusted its implementation accordingly. Following feedback from the SPORES authors, we base our implementation (Algorithm \ref{alg:spores-alg}) on the evolving-average version of SPORES from \citet{lombardi2023redundant}, while using the VMM-only initialization as in \citet{LCA_spores_2025}. 
        The update step on line 17 of Algorithm \ref{alg:spores-alg} uses equations \eqref{eq:spores-update} to update the search direction. 
        \begin{subequations}\label{eq:spores-update}
        \begin{align}
            w^{(s)}_j &= \frac{\bar{z}^{(s)}}{|\bar{z}^{(s)}-z^{*,(s)}|}, \qquad j>1 \\
            w^{(s)}_1 & =  \frac{z^{*,(s)}}{z^{max}}
        \end{align}
        Following \citet{lombardi2023redundant}, and input from the SPORES authors we weight the combination of the VMM and evolving average terms using $\alpha=1$, and $\beta=0.5$. These weights ensure that both terms are of similar magnitude, and that the VMM term does not dominate the contribution from the evolving average term. 
        \end{subequations}
    \end{description}
     \begin{algorithm}
\caption[SPORES]{SPORES implementation\\
Requires weighting parameters $\alpha$ and $\beta$, which we assign values of $1$ and $0.5$.}
\label{alg:spores-alg}
\begin{algorithmic}[1]
    \State $i \gets 1$
    \While{$i \le n_z$} \Comment{Begin each sequence with a VMM weight}
        \State $w_0^{(i)} \gets \mathbf{0}^{n_z}$
        \State $w^{(i)}_{0,i} \gets -1$ 
        \State $w_0^{(i+n_z)} \gets \mathbf{0}^{n_z}$
        \State $w^{(i+n_z)}_{0,i} \gets 1$
        \State $i \gets i + 1$
    \EndWhile
    \State $k \gets 1$ \Comment{Counter for MGA iterations}
    \While{$k \le m$}
        \State $s \gets 1 + (k \bmod 2n_z)$ \Comment{Index sequence by $s$}
        \State $j \gets \left\lfloor 1 + \frac{k}{2n_z} \right\rfloor$ \Comment{Floored division}
        \If{$j == 1$}
            \State $z^{*(s)} \gets \texttt{mga}(w^{(s)}_j)$ \Comment{Solve \eqref{eq:MGA}}
            \State $\bar{z}^{(s)} \gets z^{*(s)}$
        \Else
            \State $\bar{z}^{(s)}\gets \frac{ (j-1)\bar{z}^{(s)} +  z^{*(s)}}{j}$ \Comment{Update moving average}
            \State $w^{(s)}_j \gets \texttt{update}(w^{(s)}_{j-1}, z^{*(s)}, \bar{z}^{(s)})$
            \State $z^{*(s)} \gets \texttt{mga}(\beta w^{(s)}_0+\alpha w^{(s)}_j)$ \Comment{Solve \eqref{eq:MGA}}
        \EndIf
        \State $k \gets k + 1$
    \EndWhile
\end{algorithmic}
\end{algorithm}
    \subsection{Max-min reformulations}\label{subsec: max-min reformulations}
    By selecting $p=\infty$ and expanding the norm, problem     \eqref{eq: maxmin_dist_prob} can be written as the bi-level LP:
\begin{subequations}\label{eq: bi-level-pinfty}
    \begin{align}
        \min_{t,z_\mathcal{I},z_\mathcal{O}}\:\:& -t \\
        \text{s.t.}\quad     &z_\mathcal{O}\in \mathcal{O}\\
        & t, z_\mathcal{I} \in \arg\min_{t, z_\mathcal{I}} \quad  \{t\::\:t \le z_\mathcal{O} - z_\mathcal{I} \le t,\: z_\mathcal{I} \in \mathcal{I} \}
    \end{align}
\end{subequations}
    To solve this bi-level LP, we reformulate it into a single-level MILP using the KKT conditions of the lower level problem. This is a standard reformulation described in the literature \citep{bard2013practical, pineda2018efficiently}; however, for completeness, we briefly describe the reformulation below.
    
    Problem \eqref{eq: bi-level-pinfty} can be manipulated into the (simplified) standard form:
    \begin{subequations}\label{eqn: standard-bilevel}
        \begin{align}
            \min_{x_U,y_L}\:\:& -c_t^Ty_L \\
            \st\:\:& A_Ux_U\le b_U  \\
            & y_L \in \arg\min_{y_L} \{c_t^Ty_L\::\: B_Ly_L+C_Lx_U\le d_L,\: F_Ly_L=e_L\}
        \end{align}
    \end{subequations}
    where $x_U\inR^{n_{x_U}}$ are the upper-level variables, $y_L\inR^{n_{y_L}}$ are the lower-level variables, the upper-level constraints are written as $n_{I_U}$ inequality constraints with $A_U\inR^{n_{I_U}\times n_{x_U}}$ and $b_U\inR^{n_{I_U}}$, and the lower-level constraints are written as $n_{I_L}$ and $n_{E_L}$ inequality and equality constraints with $B_L\inR ^{n_{I_L}\times n_{y_L}}$, $C_L\inR^{n_{I_L}\times n_{x_U}}$, $d_L\inR^{n_{I_L}}$, $F_L\inR ^{n_{E_L}\times n_{y_L}}$, and $e_L\inR ^{n_{E_L}}$. $c_t\inR ^{n_{y_L}}$ is a vector used to select the element of $y_L$ corresponding to $t$.
    We reformulate problem \eqref{eqn: standard-bilevel} by replacing the lower level problem with its KKT to yield the single-level problem:  
    \begin{subequations}
        \begin{align}
            \min_{x_U,y_L}\:\:& -c_ty_L \\
            \st\:\:& A_Ux_U\le b_U  \\
            &B_Ly_L+C_Lx_U\le d_L\\
            &F_Ly_L=e_L\\
            & c_t +  B_L^T\lambda_L +  F_L^T\eta_L = 0 \label{eq: single-reform-kkt-stat}\\
            & \lambda_L\ge0 \label{eq: single-reform-dual-feas}\\
            & \lambda_{L,i}(B_{L,i}y_L+C_{L,i}x_U-d_{L,i})=0\qquad \forall i=1\dots n_{I_L} \label{eq: comp-con-need-reform}
        \end{align}
    \end{subequations}
    where $\lambda_L\inR^{n_{I_L}}$ and  $\eta_L\inR^{n_{E_L}}$ are the dual variables corresponding to the inequality and equality constraints of the lower level problem, equation \eqref{eq: single-reform-kkt-stat} is the stationary condition, equation \eqref{eq: single-reform-dual-feas} is the dual feasibility condition, and equation \eqref{eq: comp-con-need-reform} is the complementary constraint.

    Complementary constraints are (i) inherently combinatorial and (ii) difficult to handle directly as they fail to satisfy constraint qualifications assumed by general non-linear solvers. To avoid the second issue, we reformulate equation  \eqref{eq: comp-con-need-reform} using a MILP reformulation:
    \begin{subequations}
        \begin{align}
        \lambda_i &\le u_iM_D \qquad &&\forall i=1\dots n_{I_L}\\
        B_iy_L+C_ix_U-d_i &\le M_P(1-u_i)\qquad &&\forall i=1\dots n_{I_L}
        \end{align}
    \end{subequations}
    where $M_P$ and $M_D$ are sufficiently large constants that bound the primal and dual solutions \citep{bard2013practical}, and $u\in\{0,1\}^{n_{I_L}}$. If $u_i=0$, the corresponding inequality constraint is inactive, whereas if $u_i=1$, the inequality is active. 
    
    Using information about the problem, we can tighten this reformulation by providing bounds on the sum of $u$. From Carathéodory's theorem, a maximum of $n_z+1$ points are needed to interpolate any point in the convex hull. Additionally, a maximum of $n_z$ constraints can be active due to the infinity norm. Similarly, at least one point is needed for the interpolation, and at least one of the norm inequalities must be active. 
    
    \subsection{ORACLE pseudocode and variations}\label{subsec: oracle-pseudo}
 The basic ORACLE algorithm is given in Algorithm \ref{alg: sandwich-alg}, and described in Section \ref{sec: methodology}.
 ORACLE can be extended in various ways. Here, we briefly outline some simple, but potentially beneficial variations:
\begin{description}
    \item[Parallelization:] The version of ORACLE in Algorithm \ref{alg: sandwich-alg} is sequential due to the generation of only one trial point at each iteration. However, solution pools of a desired size can be returned when finding the distance between $\IA$ and $\OA$ in Step 2. These solutions are feasible points found during ``normal" optimization or specially selected by the optimization software (e.g., to maximise diversity in the pool). All points in the solution pool can be evaluated in parallel to check for feasibility and to generate cuts. Thus, a partial degree of parallelisation can be easily implemented.
    \item[Combination with other MGA algorithms:] All MGA approaches using the LP formulation \eqref{eq:MGA} find points on the perimeter of $\IA$, and can also generate constraints for $\OA$. Thus, any of these algorithms can be used within ORACLE to generate additional points and constraints. Then, ORACLE will sequentially check if convergence has been achieved (while introducing constraints), while the MGA algorithm is used to generate more points in parallel.
    \item[MILP heuristics] Step 2 requires solving an MILP, which can be computationally expensive. Thus, to further reduce the computational effort, one could accurately solve problem \eqref{eq: maxmin_dist_prob} only every $n$ iterations, and instead use local methods for the intermediate iterations. These local methods are only heuristics; however, they often perform well \citep{kleinert2021computing,bard2013practical} and will generate valid trial points for Step 3. 
    However, we note that as system models can be very computationally expensive \citep{brochin2024harder}, Step 2 of ORACLE may not be the computational bottleneck. %
\end{description}

 \begin{algorithm}
    \caption[ORACLE]{ORACLE: Optimization for the Rigorous Analysis of the Complete Landscape of Alternatives\\
    Figure \ref{fig:prop-alg} visually depicts the algorithm, with matching colours.}\label{alg: sandwich-alg}
    \begin{algorithmic}[1]
    \State ${\color{green}\mathcal{O}},{\color{blue}\mathcal{I}},{\color{orange} d_{\mathcal{OI}}} \gets \texttt{init\_regions}$\Comment{Step 1}
    \While{$ d_{\mathcal{OI}} > \texttt{tol}$}
        \State ${\color{orange}z_{\mathcal{O}}^*,\: d_{\mathcal{OI}}} \gets \texttt{furthest\_point}({\color{green}\mathcal{O}},{\color{blue}\mathcal{I}})$ \Comment{Step 2}
    \State ${\color{red}z_{f}^*, d_{\mathcal{OZ}}} \gets \texttt{closest\_near\_opt}({\color{orange}z_{\mathcal{O}}^*})$ \Comment{Step 3}
        \State ${\color{green}\mathcal{O}},{\color{blue}\mathcal{I}} \gets \texttt{refine}({\color{green}\mathcal{O}},{\color{blue}\mathcal{I}}, {\color{red}z_{f}^*,d_{\mathcal{OZ}}})$  \Comment{Step 4}
    \EndWhile
    \State ${\color{magenta}z_{designs}} \gets \texttt{sample}({\color{green}\mathcal{O}},{\color{blue}\mathcal{I}})$  \Comment{Step 5}
    \end{algorithmic}
    \end{algorithm}
    \subsection{Cutting the outer approximation}\label{subsec: hyperplane-from-dual}

    We describe two approaches to cut the outer approximation (if the trial point is infeasible): the first approach (\ref{subsec: sep-hyperplane}) ensures the discovery of a strictly separating hyperplane, whereas the second approach (\ref{subsec: double-cut-OA}) relies on a local approximation of the total cost, which requires an additional optimization. The rationale of the second approach is that if the trial point was infeasible \emph{only} due to the near-optimality constraint, then the second approach would generate a facet of the true near-optimal region. In general, it is unlikely that only the near-optimality constraint would not be satisfied at every iteration, and so the second approach cannot be used to guarantee convergence. 

    Our implementation of ORACLE uses both approaches whenever the trial point is infeasible. The additional optimization for the second approach increases the computational cost per iteration; however, adding two cuts can reduce the total number of iterations (and hence total computational cost) for convergence\footnote{Ideally, this increase in computational cost can be partially reduced by warm-starting the second optimization.}.   
    In general, the computational trade-off between additional effort per iteration and the total number of iterations depends on the system under consideration and the desired tolerance.

\subsubsection{A strictly separating hyperplane }\label{subsec: sep-hyperplane}
    A strictly separating hyperplane is given by the dual variables of equation \eqref{eq: define-delta-for-duals}. To show that the dual variables have this property, consider the dual  of problem \eqref{eq: closest-feasible-inf}:
    \begin{subequations}\label{eq: closest-feasible-inf-DUAL}
\begin{align}
\max_{\mu,\lambda_A,\eta_S,\eta_F,\lambda_v} \quad & \mu^T z_{\mathcal{O}}^* + \inf_{x,z} \Big( \lambda_A^T (Ax - b) + \eta_S^T (Sx - z_f) \notag \\
& \quad + \eta_F^T (Fx - d) + \lambda_v (c^T x - v^*(1+\epsilon))  - \mu^T z_f \Big) \label{eq: closest-feasible-inf-DUAL-OBJ} \\
\text{s.t.} \quad & \lambda_A \ge 0, \quad \lambda_v \ge 0, \quad \|\mu\|_{p^*} \le 1 
\end{align}\end{subequations}
where $\|\cdot\|_{p*}$ is the dual norm\footnote{The dual of the $l_1$ norm is the $l_\infty$ norm and vice versa. From  Hölder's inequality the dual of the $l_p$ is the $l_q$ norm  with $q=\frac{p}{p-1}$ for $p,q\in (1,\infty)$.} of $\|\cdot\|_{p}$, and $\mu\inR^{n_z}$ is the dual variable of \eqref{eq: define-delta-for-duals}, and $\lambda_A,\eta_S,\eta_F,\lambda_v$ are the other dual variables.

Suppose that $p$ is chosen so that problem \eqref{eq: closest-feasible-inf} is convex, and that strong duality holds. Let $\mu^*,\lambda_A^*,\eta_S^*,\eta_F^*,\lambda_v^*$ correspond to a dual optimum solution, with corresponding primal solution $z_f^*, x^*$.
If $z_{\mathcal{O}}^* \not=z_f^*$, the optimal objective of problem \eqref{eq: closest-feasible-inf} is positive. By strong duality, the optimum value of equation \eqref{eq: closest-feasible-inf-DUAL} is equal to the primal optimum value and hence also positive. If the optimal dual objective is positive, then when using the optimum dual variables:
 \begin{subequations}
          \begin{align}
          \eqref{eq: closest-feasible-inf-DUAL-OBJ} &>0\\
            \implies \mu^{*T}z_{\mathcal{O}}^*-\inf_{z_f} \mu^{*T}z_f&>0\\
            \implies \mu^{*T}z_{\mathcal{O}}^*-\mu^{*T}z_f&>0\qquad \forall \text{ primal feasible $x,z_f$}
          \end{align}
\end{subequations}
Thus, $\mu^{*}$ defines the normal of a hyperplane that strictly separates $z_{\mathcal{O}}^*$ and the near-optimal set.

\subsubsection{Cutting the outer approximation using the total system cost}\label{subsec: double-cut-OA}

After a feasible point $z_f^*$ has been found, one can generate a local approximation of the optimal total system cost and add this constraint to the outer approximation. The benefit of introducing this second cut is that the convergence of the outer approximation may be accelerated, and hence the total number of iterations may be reduced. Consider the parametric LP:
 \begin{subequations}\label{eq: total-cost-cut-prob}
          \begin{align}
             v_f^*(z) = \min_{x}\:\: & c^T x\\
            \st\:\: &Sx = z \label{eq: constr-cap}\\
            & Ax \leq b\\
            &  Fx = d
          \end{align}
\end{subequations}
which finds the cost-optimal solution, subject to constraint \eqref{eq: constr-cap}, which fixes the exploratory variables to some vector $z$, i.e., the remaining decision variables are cost-optimal. 

Importantly, the value function $v_f^*$ is a convex, piecewise linear function of the parameter $z$. 
The dual variable of constraint \eqref{eq: constr-cap}, $\lambda\inR^{n_z}$, can be used to construct a linear approximation of the value function, parametrized by $z$. Consider this linearization for the point $z=z_f^*$:
\begin{equation}\label{eq: linear-value-approx}
     \hat{v}_f^*(z) \approx  v_f^*(z_f^*)+ \lambda^T(z_f^*)(z-z_f^*) 
\end{equation}
As $\lambda(z_f^*)$ defines the derivative of $v_f^*$ with respect to $z_f^*$,  convexity implies that $ \hat{v}_f^*$ is a global under-estimator of $v_f^*$. Therefore the constraint:
\begin{equation}
    v_f^*(z_f) + \lambda(z_f^*)^T(z-z_f^*)  \le  v^*(1+\epsilon)
\end{equation}
is a supporting hyperplane of the near-optimal region. Additionally, if $z_{\mathcal{O}}^*$ is infeasible only due to the near-optimality constraint, the normal of this constraint is a normal of a facet of the near-optimal region.

If the system is modelled such that any $z$ (within its upper and lower bounds) is feasible if the near-optimality constraint is neglected, $z_f^*$ will be the cost-optimal solution, as otherwise more of the near-optimal budget could be ``spent" to move $z_f^*$ closer to $z_{\mathcal{O}}^*$, i.e. the decisions $x$ will be cost-optimal. Thus, the primal solution of problem \eqref{eq: closest-feasible-inf} could be used to construct the primal and dual solution of problem \eqref{eq: total-cost-cut-prob} without requiring the optimization problem to be solved. In general, the system will not necessarily be modelled in this way, and problem \eqref{eq: total-cost-cut-prob} will have to be solved. However, the optimizer can be warm-started by using the solution of problem \eqref{eq: closest-feasible-inf}.

\subsection{Comparing convergence by volume ratio}\label{subsec: compare-volume-ratio}
\begin{figure}
    \centering\includegraphics[width=.99\linewidth]{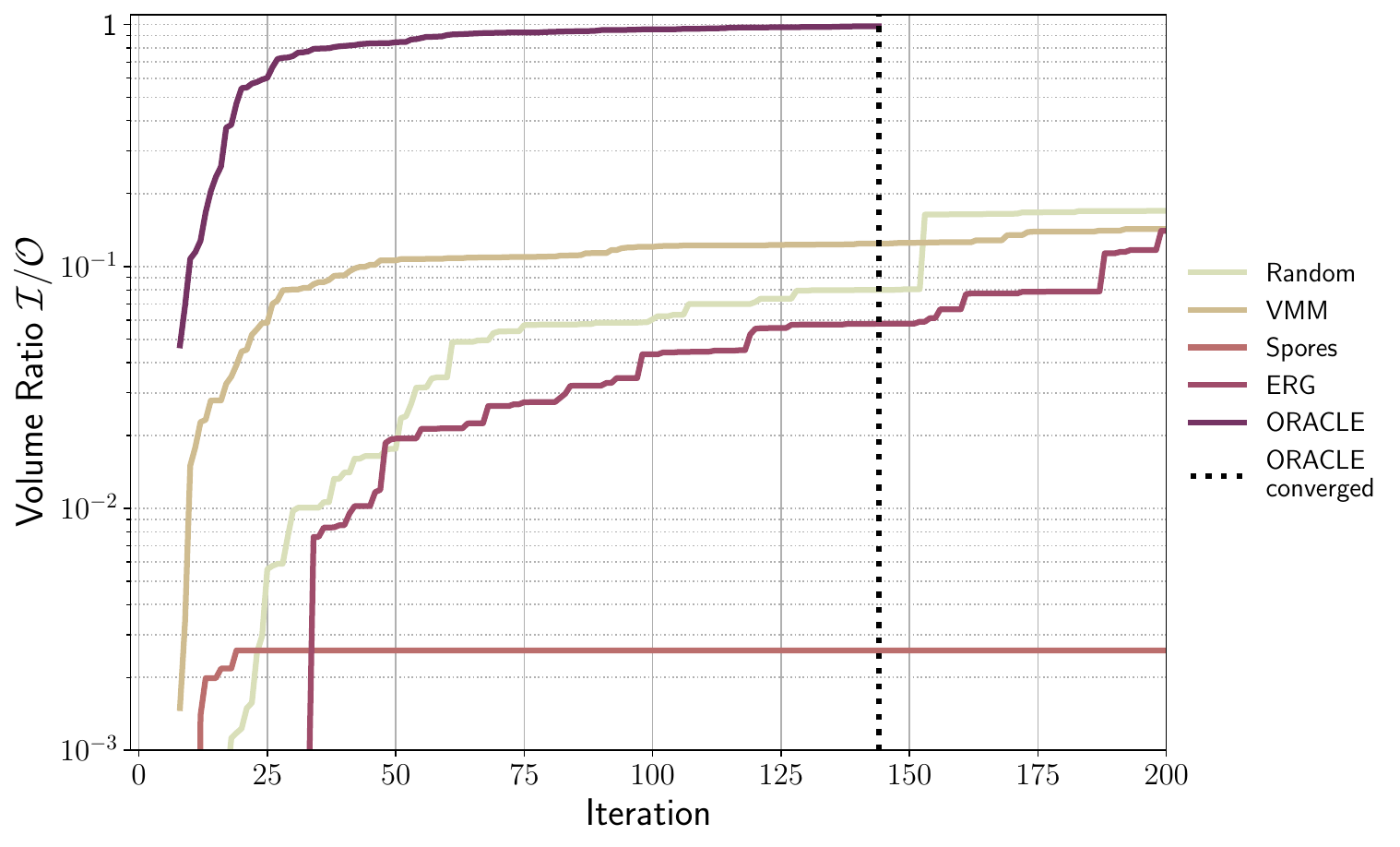}
    \caption{Ratio of the volume of the inner and outer approximations of state-of-the-art MGA methods and ORACLE. ORACLE terminates after 144 iterations, having converged to within 0.1 GW, which corresponds to a volume ratio of 0.98.}
    \label{fig:compare-volumes-small-tech-list-ratio}
\end{figure}

\citet{schricker2023gotta} propose monitoring the ratio of volumes, as illustrated in Figure~\ref{fig:compare-volumes-small-tech-list-ratio}. This approach yields the same ranking as our proposed distance-based metric; however, the volume ratio is generally less interpretable. The same ranking is achieved as both \citet{schricker2023gotta} and our metric use inner and outer approximations, which enables the tracking of convergence. While volume-based metrics are limited to low-dimensional cases ($n_z\lesssim8$) \citep{barber2013qhull}, our proposed metric is applicable to higher-dimensional systems. %

\subsection{Comparing convergence to the near-optimal region}\label{subsec: compare-distance-true}

\begin{figure}
    \centering
    \includegraphics[width=0.99\linewidth]{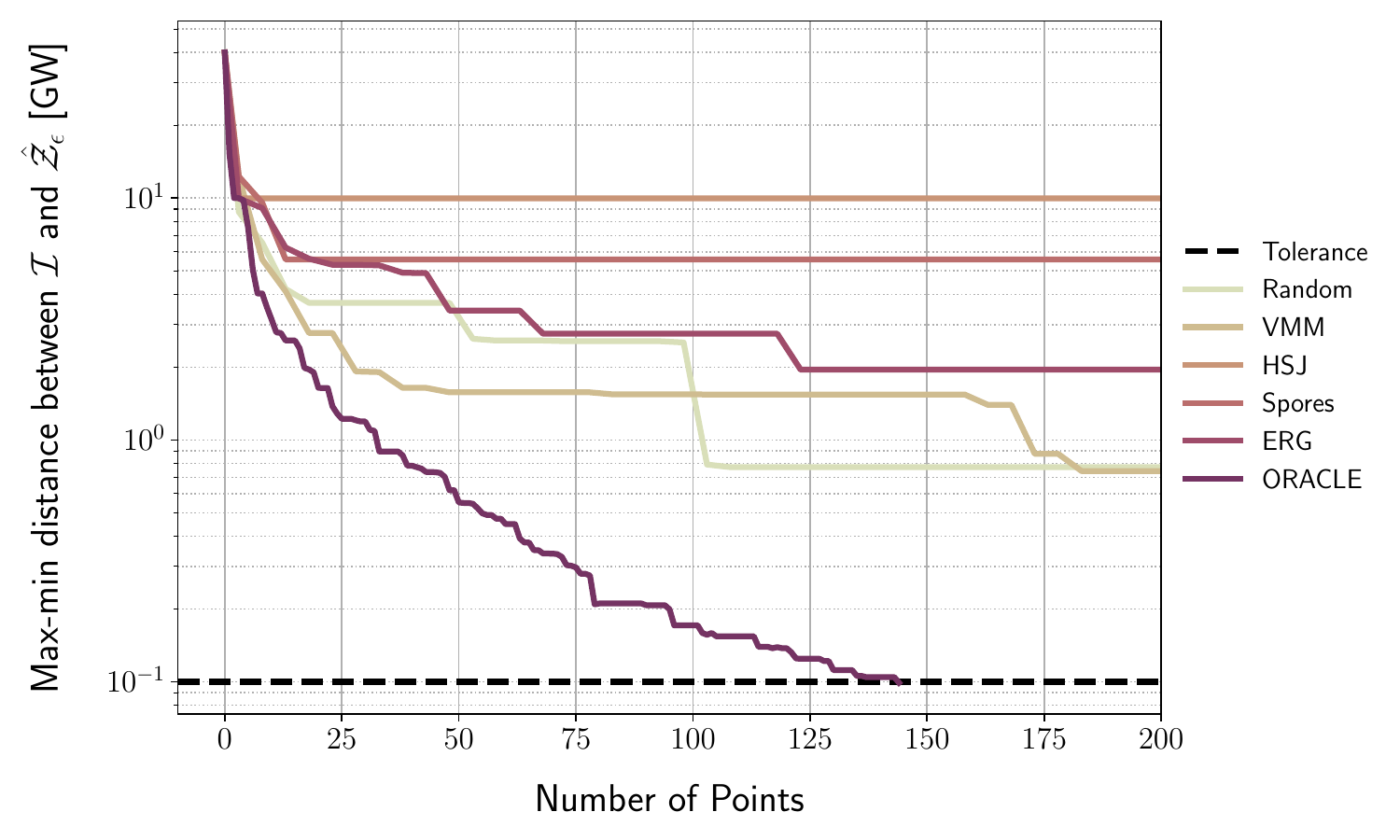}
    \caption{Approximate distance between the true near-optimal region and the inner approximations of state-of-the-art MGA methods and the ORACLE algorithm.}
    \label{fig: dist-conv-true-region}
\end{figure}

Figures \ref{fig:distance-oracle-small-tech-list}, \ref{fig:compare-volumes-small-tech-list}, and \ref{fig:compare-volumes-small-tech-list-ratio} plot the convergence of the compared method's inner and outer approximations at iteration $k$, using the data generated by that method up to iteration $k$. As the final outer-approximation of ORACLE has converged to within 0.1 GW of the true near-optimal region, we can use it as an approximation of the true region, $\hat{\mathcal{Z}}_\epsilon$. Thus, in a post-processing step, we can measure the distance of any method's inner approximation to the true region (Figure \ref{fig: dist-conv-true-region}).

Despite the inner-approximations of ORACLE, ERG, VMM and Random appearing to have converged to nearly the same value in Figure \ref{fig:compare-volumes-small-tech-list}, Figure \ref{fig: dist-conv-true-region} shows that the inner-approximations of these methods have not converged to within the specified tolerance. By 200 iterations, VMM and Random have converged to within around 0.8 GW, and ERG to 2 GW. These distances are smaller than the distances that the methods can report using only the data they generate (Figure \ref{fig:distance-oracle-small-tech-list}). This is consistent with the finding that the outer approximations of Random, VMM, and ERG are not tight (see section \ref{subsec: comparison-of-volume-metrics}). 

\subsection{Computational time of ORACLE before speed up}\label{subsec: comp-time}

\begin{figure}
    \centering\includegraphics[width=.99\linewidth]{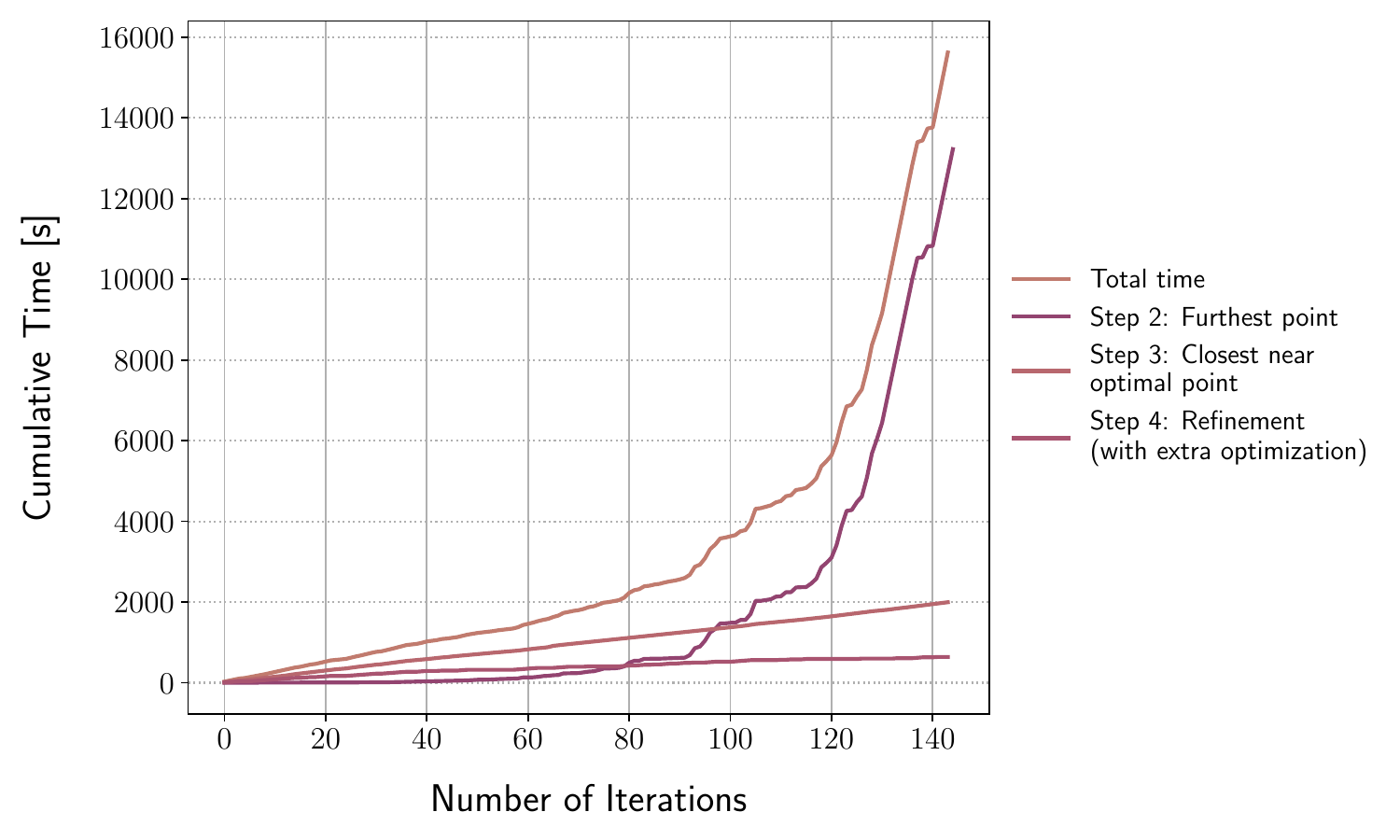}
    \caption{Cumulative computational time of the steps of ORACLE. The time of Steps 3 and 4 depends on the model used, while the time for Step 2 is primarily determined by the desired convergence tolerance and number of exploratory variables.}
    \label{fig:cum-solve-time}
\end{figure}

Figure \ref{fig:cum-solve-time} shows the contribution of each step of ORACLE to the cumulative solve time. In total, ORACLE takes under 16 000 seconds (4 hours, 26 minutes) to converge. After approximately 100 iterations, Step 2, solving the MILP, becomes the most computationally expensive contribution. However, the computational time of Steps 3 and 4 depends on the energy system model. Step 4 is computationally cheaper than Step 3 as the optimization is warm-started. The computational cost of these steps is relatively constant, as these problems do not scale with an increasing number of iterations. Additionally, note that we specified a relatively tight tolerance of 0.1 GW for the convergence of ORACLE. If a tolerance of 0.5 GW had been specified, then ORACLE would have converged after around 55 iterations, and Step 3 would have consistently been the most computationally expensive step.

Step 2 becomes the most computationally expensive step relatively early as we are using the EnergyScope model, which was specifically designed for rapid evaluation (under 60 seconds per run) \citep{limpens2019energyscope}. Indeed, this fast run-time is specifically why the model was chosen, as we wished to compare various methods to each other. 
In contrast, most energy system models require significantly longer evaluation times due to more detailed representations; for example, PyPSA-Eur  \citep{PyPSAEur} takes around 4 hours (14 400 seconds) per optimization. 
The computational cost of Step 2 is primarily governed by the number of iterations (which is determined by the convergence tolerance) and the number of exploratory variables. The complexity of the model does not determine these factors. Therefore, in a study with a model similar to PyPSA-Eur, it appears reasonable to expect that the total computational burden of Step 2 will be negligible.  

\end{document}